\documentclass[10pt,a4paper]{article}
\begin{document}
\leftline{\bf Old and new algorithms for $\pi$}	
\thispagestyle{empty}

This letter concerns Semjon Adlaj's article \emph{An eloquent formula for
the perimeter of an ellipse} [Notices 59, 8 (Sept.\ 2012), 1094--1099].
In his comments on the ``(so-called) Brent-Salamin algorithm'' for
computing $\pi$, Prof.\ Adlaj misses some important points.

First, both Brent and Salamin acknowledged their debt to Gauss and Legendre.
That the names ``Brent-Salamin'' or ``Salamin-Brent''
are widely used is probably due to the ambiguity of calling something new
after Gauss and Legendre, e.g.\ a Google search for ``Gauss-Legendre''
gives many hits on Gauss-Legendre quadrature.

Second, although Euler discovered the special case of Legendre's relation
that is used in the simplest Brent-Salamin algorithm ($k = k' =
1/\sqrt{2}$), the more general form of Legendre's relation is needed for the
members of the family of algorithms that arise from choosing $k \ne k'$. 
Since Legendre's relation is not attributed to Euler, it would be
uninformative to use the name ``Gauss-Euler'' as Prof.\ Adlaj suggests
[footnote 4]. A Google search for ``Gauss-Euler'' gives even more hits than
one for ``Gauss-Legendre'', but they are almost all irrelevant.

Third, and more important, none of those three great mathematicians of the
past would have appreciated the significance of such an algorithm, because
they lived in the days before electronic computers and fast algorithms, such
as the Sch\"onhage-Strassen algorithm, for multiplication of large integers. 
Without such technology and modern algorithms, the Brent-Salamin algorithm
is a relatively poor algorithm for computing $\pi$~-- algorithms based on
the Maclaurin series for $\arctan(1/n)$, such as Machin's
$\pi/4 = 4\arctan(1/5) - \arctan(1/239)$, are far superior (even today,
they are competitive if combined with binary splitting and fast
multiplication algorithms).  Indeed, on reading Gauss's unpublished notebook
entry of May 1809, it seems probable that he did not regard his discovery as
an algorithm for computing $\pi$, since $\pi$ only appears in the
denominator of the right-hand side of the crucial equation. More likely
Gauss regarded this equation as an interesting identity involving elliptic
integrals, only incidentally involving the known constant $\pi$. [The
relevant notebook entry is reproduced on page 99 of the book \emph{Pi:
Algorithmen, Computer, Arithmetik} by Arndt and Haenel.]

Finally, perhaps this emphasis on the computation of a single constant is
unwarranted. Brent's 1975 and 1976 papers, not referenced by Prof.\ Adlaj,
showed that \emph{all} elementary functions can be evaluated to
given accuracy just as fast as $\pi$, up to a constant factor, by using the
arithmetic-geometric mean.  This of course includes the computation of an
infinite set of constants such as
$e^{\pi}$ and $\pi/e$. No doubt this fact would have been of more
interest to Euler, Legendre and Gauss than yet another formula or algorithm
for $\pi$.
\smallskip

\rightline{--- Richard Brent}
\rightline{Australian National University}
\rightline{adlaj@rpbrent.com}
\end{document}